\documentclass{amsart}
\usepackage{amssymb}
\usepackage{amsmath}

 \usepackage [latin1]{inputenc}

\newtheorem{theorem}{Theorem}

\newtheorem{example}[theorem]{Example}

\usepackage{graphicx}

\def\eps{\varepsilon}
\graphicspath{%
    {converted_graphics/}
    {/}
}
\begin{document}
\title[An Extreme Metastable Model for Quantum Measurement]
{  An Extreme Metastable Model for Quantum Measurement}

\author[A. Boyarsky]{Abraham Boyarsky }
\address[A. Boyarsky]{Department of Mathematics and Statistics, Concordia
University, 1455 de Maisonneuve Blvd. West, Montreal, Quebec H3G 1M8, Canada}
\email[A. Boyarsky]{abraham.boyarsky@concordia.ca}

\author[P. G\'ora]{Pawe\l\ G\'ora }
\address[P. G\'ora]{Department of Mathematics and Statistics, Concordia
University, 1455 de Maisonneuve Blvd. West, Montreal, Quebec H3G 1M8, Canada}
\email[P. G\'ora]{pawel.gora@concordia.ca}

\author[Z. Li]{Zhenyang Li }
\address[Z. Li]{Department of Mathematics, Honghe University, Mengzi, Yunnan 661100, China}
\email[Z. Li]{zhenyangemail@gmail.com}

\thanks{The research of the authors was supported by NSERC grants. The research of Z. Li is also supported
by NNSF of China (No. 11161020 and No. 11361023)}

\date{\today }
\keywords{Quantum measurement problem, metastable dynamical system, density functions, extreme points}

\begin{abstract}  
Quantum experiments are observed as probability density functions. We often encounter multi-peaked densities which we model in this paper by a metastable dynamical system. The dynamics can be regarded in a thought experiment where a mouse is in either one of two disjoint traps which possess little outlets. The mouse spends most of its time in one or the other trap, and once in awhile makes its way to the other trap. Hence, a bi-peaked density function. When the experiment is observed - such as a light shone on the traps - the mouse stays in one trap. This measurement process results in a single peaked density that, we prove, can be modelled by the metastable process selecting an extreme density function, which is single peaked.
\end{abstract}

\maketitle

\textbf{Highlights:}

1. Many quantum mechanical processes display a multi-peaked density function.

2. We model the underlying dynamics by a metastable dynamical system based on a piecewise expanding map.

3. A measurement on the process results in a single peaked density function, which we prove is an extreme metastable system.

\bigskip
\bigskip

{Department of Mathematics and Statistics, Concordia University, 1455 de
Maisonneuve Blvd. West, Montreal, Quebec H3G 1M8, Canada}

and

Department of Mathematics, Honghe University, Mengzi, Yunnan 661100, China

\smallskip

E-mails: {abraham.boyarsky@concordia.ca}, {pawel.gora@concordia.ca}, {%
zhenyangemail@gmail.com}.

\section{ Introduction}\label{section_intro}

The behavior of quantum systems is governed by the Schrodinger`s Wave
Equation. The complex-valued solutions of these equations are called wave
functions. When squared, the wave equation describes the probability of a
quantum particle being found in various locations of the state space.
Quantum particles have no fixed location. Their existence is spread
probabilistically across the domain of the wave function. However, when a
macro measurement is made there is an instantaneous and dramatic change: all
random effects vanish and real properties immediately appear. Quantum theory
does not explain why and what happens as probability turns into certainty.
This transition is called collapse of the wavefunction and the mysterious
effect is referred to as the measurement problem.

\begin{figure}[h] 
  \centering
  \includegraphics[bb=0 0 400 400,width=2.22in,height=2.22in,keepaspectratio]{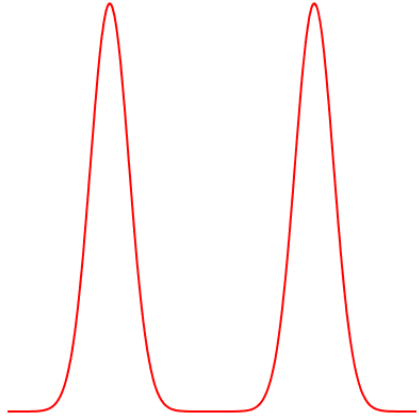}
  \caption{Probability density function with two narrow peaks}
  \label{fig:2peaks}
\end{figure}

A common way of portraying the measurement problem is to display the
probability density function (pdf) of finding a quantum particle. In the
quantum case it is a pdf $f\ $ with two narrow peaks \cite[page 31]{timfo}. When a
measurement is made one or other of the two peaks vanishes and the quantum
particle is localized to the remaining single peak region of the state space.
There are various explanations for this effect: for example, decoherence and
continuous spontaneous localization, but neither are completely satisfactory.

\begin{figure}[h] 
  \centering
  \includegraphics[bb=0 0 400 400,width=2.22in,height=2.22in,keepaspectratio]{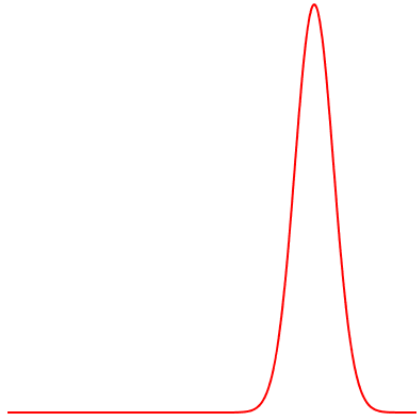}
  \caption{Probability density function $f_0$ with one narrow peak}
  \label{fig:1peaks}
\end{figure}

In this note we will model the measurement problem by using a deterministic
metastable system. First, we refer to the following thought experiment: a
mouse is moving about in a set $A$ which has a small escape hatch. The mouse
spends most of its time in $A$, but now and again it escapes and makes its
way to a set $B$ where it behaves in the same way as in $A$. To locate the
mouse a light is shone on the setup. The process of shining light is a macro
measurement; it frightens the mouse and it no longer ventures out of the set
it is in at the time of measurement.

In section 2 we describe a general discrete time metastable system that
models the above mouse trap experiment. At the outset we assume that a map $T
$ satisfies the conditions of \cite{CBP} that are sufficient for for its
deterministic perturbation $T_{\varepsilon }$ to be metastable  as well as
generating the pdf $f$ \ described above. Note that our starting point is
not the wave function $\psi $ but rather the two peaked pdf $f$ which stems
from a wave function $\psi $. We do not attempt to relate $T$ to $\psi $,
although both generate $f$ but in completely different ways. 

Let $A$ and $B$ be two disjoint sets of the state space. $T$ is assumed to
have an absolutely continuous invariant measure (acim) $\mu $ (not ergodic)
on $A\cup B.$ Now we assume that there are two mutually singular ergodic
acims $\mu _{A}$ and $\mu _{B},\ $ on $A$
and $B,$ respectively, where $\mu _{A}$ and $\mu _{B}$ have probability
density functions (pdf) $f_{A}$ and $f_{B}$. Rather than refer back to the
wave functions that generate these pdfs, we construct a piecewise smooth
expanding map $T$ that has invariant densities $f_{A}$ and $f_{B}.$To model
the motion of a quantum particle which, like the mouse in the thought
experiment, can be found in either $A$ or $B,$ we define a metastable map $%
T_{\varepsilon }$ which is an $\varepsilon $ determinsitic perturbation of $%
T.\ T_{\varepsilon \text{ }}$causes $A$ and $B$ to forfeit their invariance
and creates a single acim $\mu _{\epsilon }$ on $A\cup B.$ The measure $\mu
_{\epsilon }$ has a pdf $f_{\varepsilon }$ which converges to a convex
combination of $f_{A}$ and $f_{B}$ as $\varepsilon \rightarrow 0$.  The set
of densities $\mathcal{D}$ that can be attained by the metastable system
depends on the limit of the ratio of hole sizes in $A$ and $B$ and includes
all possible convex combinations of $f_{A}$ and $f_{B}.$ We propose that a
macro measurement of a quantum system is an extreme event much as shining
light in the mouse experiment. In Section 3 this leads us to search for the 
extreme points of 
$\mathcal{D}$, which are $f_{A}$ or $f_{B}$.

\bigskip

\section{ Metastable Dynamical System}\label{section_meta}

We assume the  map  $T:A\cup B \circlearrowleft$  is a piecewise $C^2$,
uniformly expanding on a partition set 
\begin{equation*}
\mathcal{C}_0=\{c_{0,0}<c_{1,0}<\cdots< c_{n,0}\}
\end{equation*}
which is called a critical set. $T$ has two invariant sets, $A$ and $B$. 
The point $b=A\cap B$ is called a boundary point and  the points in $H_{0}:=T^{-1}(\{b\})%
\setminus \{b\}$ are called infinitesimal holes. The metastable perturbed
system $T_\epsilon: A\cup B \circlearrowleft$ of $T$, which is also
piecewise $C^2$, has partition set 
\begin{equation*}
\mathcal{C}_\epsilon=\{c_{0,\varepsilon}<c_{1,\epsilon}<\cdots<
c_{n,\epsilon}\}.
\end{equation*}
The following properties are assumed:

\begin{itemize}
\item (I1) Unique acim on the initial invariant set:

$T|_A$ ($T|_B$) has only one acim $\mu_A$ ($\mu_B$) whose density
 is denoted by $f_A$ ($f_B$).

\item (I2) No return of the critical set to the infinitesimal holes:

for every $k > 0$ , $T^k(\mathcal{C}_0)\cap H_0=\emptyset $. 

\item (I3) Positive acims at infinitesimal holes:

$f_{A}$ is positive at each of the points in $H_{0}\cap A$ , and $f_{B}$ is
positive at each of the points in $H_{0}\cap B$. 

\item (I4) Restriction on periodic critical points. One of the following
holds:

(I4a) $\inf_{A\cup B\setminus\{\mathcal{C}_0\}} |T^{\prime }(x)|>2$;

(I4b) $T$ has no periodic critical points, except possibly a fixed point at
0 or 1. 

\item (P1) Unique acim:

for each $\epsilon> 0$, $T_\epsilon$ has only one acim $\mu_\epsilon$  with
density $f_\epsilon$. 

\item (P2) Boundary condition:

the boundary point does not move, and no holes are created near the
boundary. To be precise:

(P2a) if $b \not\in \mathcal{C}_0$, then necessarily $T(b)=b$, and we assume
further that for all $\epsilon> 0$, $T_\epsilon(b)=b$;

(P2b) if $b\in \mathcal{C}_0$, we assume that $T(b_-)< b < T(b_+)$ and also
that $b\in \mathcal{C}_\epsilon$ for all $\epsilon>0$.
\end{itemize}

Then, we have the following properties:

$T_\epsilon$ has only one ergodic acim $\mu_\epsilon$ on $A\cup B$ with
pdf $f_\epsilon$. We have two holes 
\begin{equation*}
H_{A,\epsilon}:=A\cap T_\epsilon^{-1}(B) \text{ and } H_{B,\epsilon}:=B\cap
T_\epsilon^{-1}(A) 
\end{equation*}
in $A$ and $B$, respectively, through which the orbit of $T_\epsilon$
escapes from one set to the other set. Once an orbit enters a hole, it leaves
one of the invariant sets for $T$ and continues in the other. As $%
\epsilon\rightarrow 0$, the holes converge to the place from which they
arise, which are called infinitesimal holes. By the conditions above, both $%
\mu_A(H_{A,\epsilon})$ and $\mu_B(H_{B,\epsilon})$ converge to 0 as $%
\epsilon\to 0$.

\begin{theorem}
\cite{CBP} Consider the family of perturbations \ $\left\{ T_{\varepsilon
}:\varepsilon \gtrdot 0\right\} $ of\ $T$ under the assumptions (I1-I4) and
(P1-P2) stated in Section 2.1 of \cite{CBP}. Suppose that
\begin{equation*}
l.h.r=\lim_{\varepsilon \rightarrow 0}\frac{\mu _{B}(H_{B,\epsilon })}{\mu
_{A}(H_{A,\epsilon })}
\end{equation*}
exists. Then
\begin{equation*}
f_{\varepsilon }\rightarrow \alpha f_{A}+(1-\alpha )f_{B},
\end{equation*}
as $\varepsilon \rightarrow 0$, where $\frac{\alpha }{1-\alpha }=l.h.r.$
\end{theorem}

\section{ Extreme Metastable Systems}\label{section_extreme}

The extreme values of all $l.h.r$ are $0$ and $\infty .$ Hence the extreme
points of the set of densities $\left\{ f_{\varepsilon }:\varepsilon \gtrdot
0\right\} $ occur if

1) $l.h.r.=0,$ then $f_{\varepsilon }\rightarrow f_{B}$ in $L^{1}$ as $%
\varepsilon \rightarrow 0$ or if $l.h.r.=\infty ,$ then $\ f_{\varepsilon
}\rightarrow f_{A}$ in $L^{1}$as $\varepsilon \rightarrow 0$. Thus, $f_{A}$
and $f_{B}$ are the extreme points of $\left\{ f_{\varepsilon }:\varepsilon
\gtrdot 0\right\} .$and correspond to what happens once a measurement of a
quantum system is made.

We have modeled the act of measurement by an extreme state for the
metastable system and have shown that this forces the system to display
either $f_{A}$ od $f_{B}$. The macro measurement effectively shuts the escape holes.
The measurement process has caused a collapse of the wave function: a two
peaked density function has become a one peaked density, depending on which
almost invariant set the mouse was in at the time the measurement was taken.

\section{ Construction of the map $T$}

In this section we construct an example of a metastable map satisfying the conditions of
Section \ref{section_meta} and preserving a density which is a rough approximation of the 
density $f$ of Figure \ref{fig:2peaks}. A specific example is presented at the end of the section.

We start by constructing a map $g:B\to B$, $B=[0,1]$, preserving a piecewise constant approximation of 
an arbitrary pdf $f_0$. Let $\mathcal P=\{I_i\}_{i=1}^N$ be a partition of $B$ into $N$ equal subintervals.
Let us define the  vector $v=[v_1,v_2,\dots, v_N]$, where  $v_i= \int_{I_i} f_0 dm $. Then, the function
$$f_{0a}(x)=\sum_{i=1}^N N\cdot v_i\cdot\chi_{I_i}(x), $$
is a piecewise constant approximation of $f_0$.

Let us choose $0<\beta_i<1$ and set $\alpha_i=v_i(1-\beta_i)$, $i=1,\dots,N$. Let
$$M=\left[\begin{matrix}  \beta_1&0&\dots&0&0\\
                         0&\beta_2&0&\dots&0\\
                         \vdots&\vdots&\ddots&\vdots&\vdots\\
                             0&0&\dots&\beta_{N-1}&0\\
                         0&0&\dots&0&\beta_N
             \end{matrix}\right]+\frac 1{\sum_{i=1}^N \alpha_i}\left[\begin{matrix} 1- \beta_1\\
                          1- \beta_2\\
                         \vdots&\\
                             1-\beta_{N-1}\\
                         1-\beta_N
             \end{matrix}\right]\left[\alpha_1,\alpha_2,\dots,\alpha_N\right].$$
The matrix $M$ is row-stochastic and the vector $v$ is its left invariant vector. This fact was used in \cite{Rog} to
introduce another method of constructing a piecewise linear map $g$ preserving the pdf $f_{0a}$.
We use the matrix $M=\left[m_{i,j}\right]_{1\le i,j\le N}$ to define $g$ which transforms linearly
the subinterval of length $m_{ij}/N$ of $I_i$ increasingly onto the interval $I_j$ if $i$ is odd and
 decreasingly onto the interval $I_{N-j+1}$ if $i$ is even. For more details and a general theory of such semi-Markov maps see \cite{GB1}. We make the construction slightly more complicated than usual to allow the map
$g$ to be continuous. Map $g$ is piecewise expanding and piecewise onto so pdf $f_{0a}$ is the unique $g$-invariant pdf
(for the general theory of piecewise expanding maps see \cite{BG} or \cite{LM}).

Now we can define the map $T:[-1,1]\to [-1,1]$ by
$$T(x)=\begin{cases} g(x) & ,\ \text{for}\  x\in [0,1];\\
                     -g(-x) & ,\ \text{for}\  x\in [-1,0].
\end{cases}
$$
The map $T$ has two ergodic components and preserves two invariant pdfs with disjoint supports:
$f_B(x)=f_{0a}(x)$ on $B=[0,1]$ and $f_A(x)=f_{0a}(-x)$ on $A=[-1,0]$. The point $\{0\}=A\cap B$ is the fixed point of $T$.

We now open small gates around $-1$ and $1$: the new map $T_\eps$ is equal to $T$ on $[-1+\eps,1-\eps]$ and
we have $T_\eps([-1,-1+\eps])=[0,a(\eps)]$ and $T_\eps([1-\eps,1])=[-a(\eps),0]$.
The map $T_\eps$ is metastable and preserves a unique density $f_\eps$ which is very close to 
$\frac 12 f_A+\frac 12 f_B$ since $T_\eps$ is symmetric and both $T$ and $T_\eps$ are piecewise expanding
(\cite{BG}).

\begin{figure}[h] 
  \centering
  \includegraphics[bb=0 0 400 400,width=2.22in,height=2.22in,keepaspectratio]{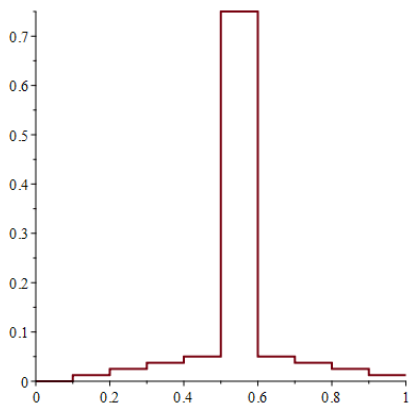}
  \caption{Density $f_B=f_{0a}$ invariant for map $g$}
  \label{fig:density_f0a}
\end{figure}

\bigskip

\begin{figure}[h] 
  \centering
  \includegraphics[bb=0 0 587 595,width=2.22in,height=2.25in,keepaspectratio]{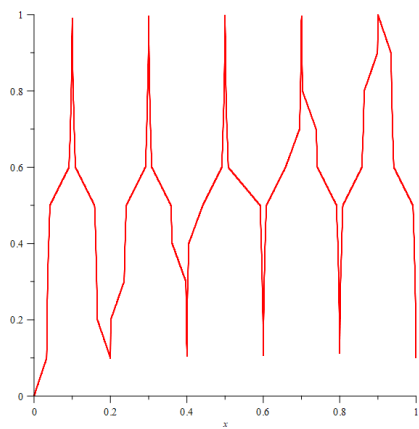}
  \caption{Semi-Markov map $g$ preserves density $f_B$}
  \label{fig:Map_g_10}
\end{figure}

\begin{example}
\end{example}
For a specific example we fix $N=10$. We consider density $f_B=f_{0,a}$ corresponding to the vector
$v=\frac 1{80}[0,1,2,3,4,60,4,3,2,1]$ (Figure \ref{fig:density_f0a}).
The semi-Markov map $g$ constructed as above is shown in Figure \ref{fig:Map_g_10}.

\newpage

\begin{figure}[h] 
  \centering
  \includegraphics[bb=0 0 958 887,width=2.22in,height=2.05in,keepaspectratio]{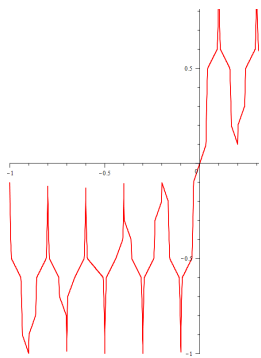}
  \caption{Map $T$ preserves two densities: $f_B$ and symmetric $f_A$}
  \label{fig:Map_bigG_10}
\end{figure}

\begin{figure}[h] 
  \centering
  \includegraphics[bb=0 0 956 863,width=2.22in,height=2in,keepaspectratio]{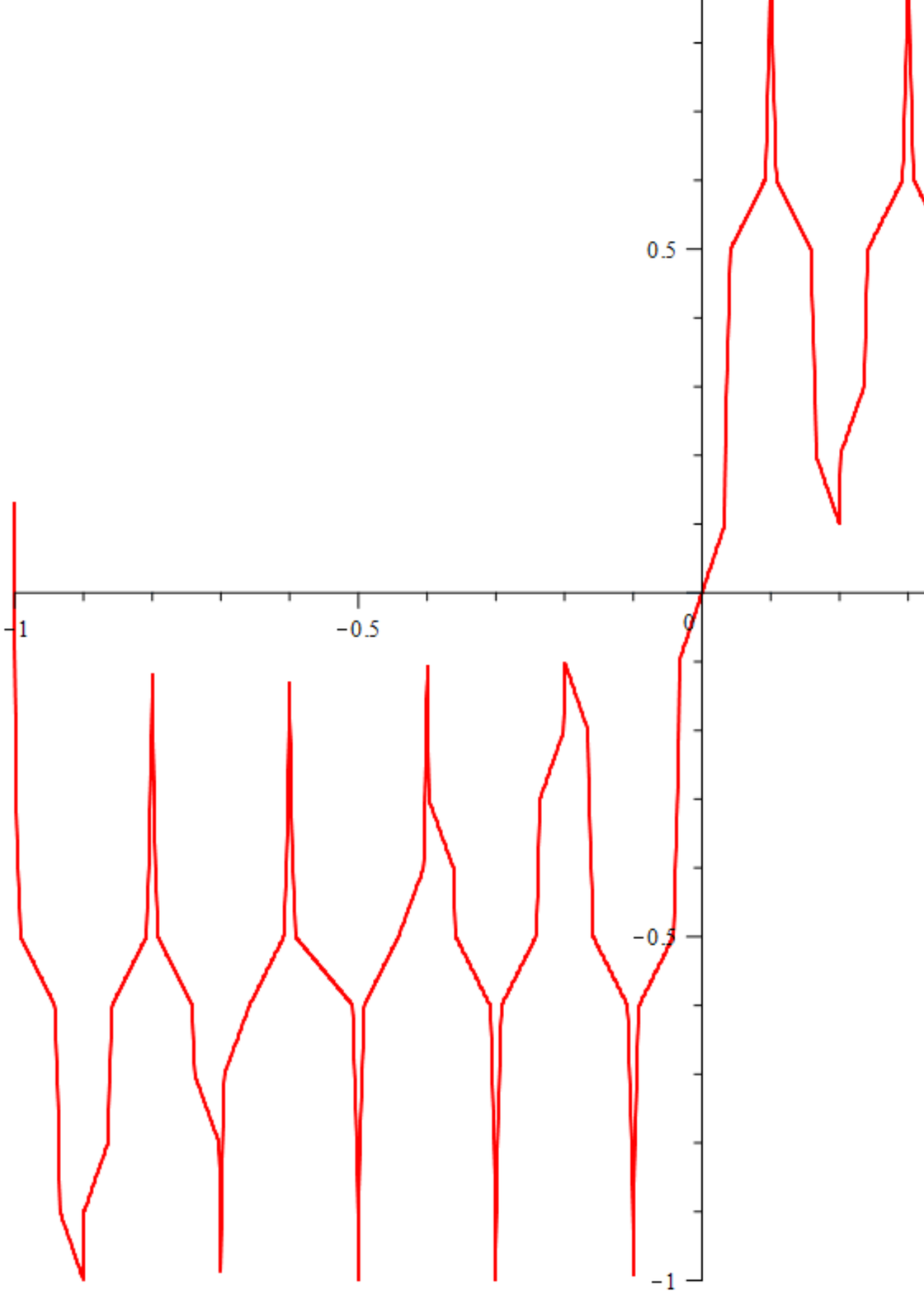}
  \caption{Metastable map $T_\eps$ preserves unique density\hskip 3 cm  $f_\eps \sim \frac 12 f_A+\frac 12 f_B$}
  \label{fig:Map_metastable_10}
\end{figure}

The maps $T$ and $T\eps$ are shown in Figures \ref{fig:Map_bigG_10} and \ref{fig:Map_metastable_10}.

\newpage

\begin{figure}[h] 
  \centering
  \includegraphics[bb=0 0 400 400,width=2.22in,height=2.22in,keepaspectratio]{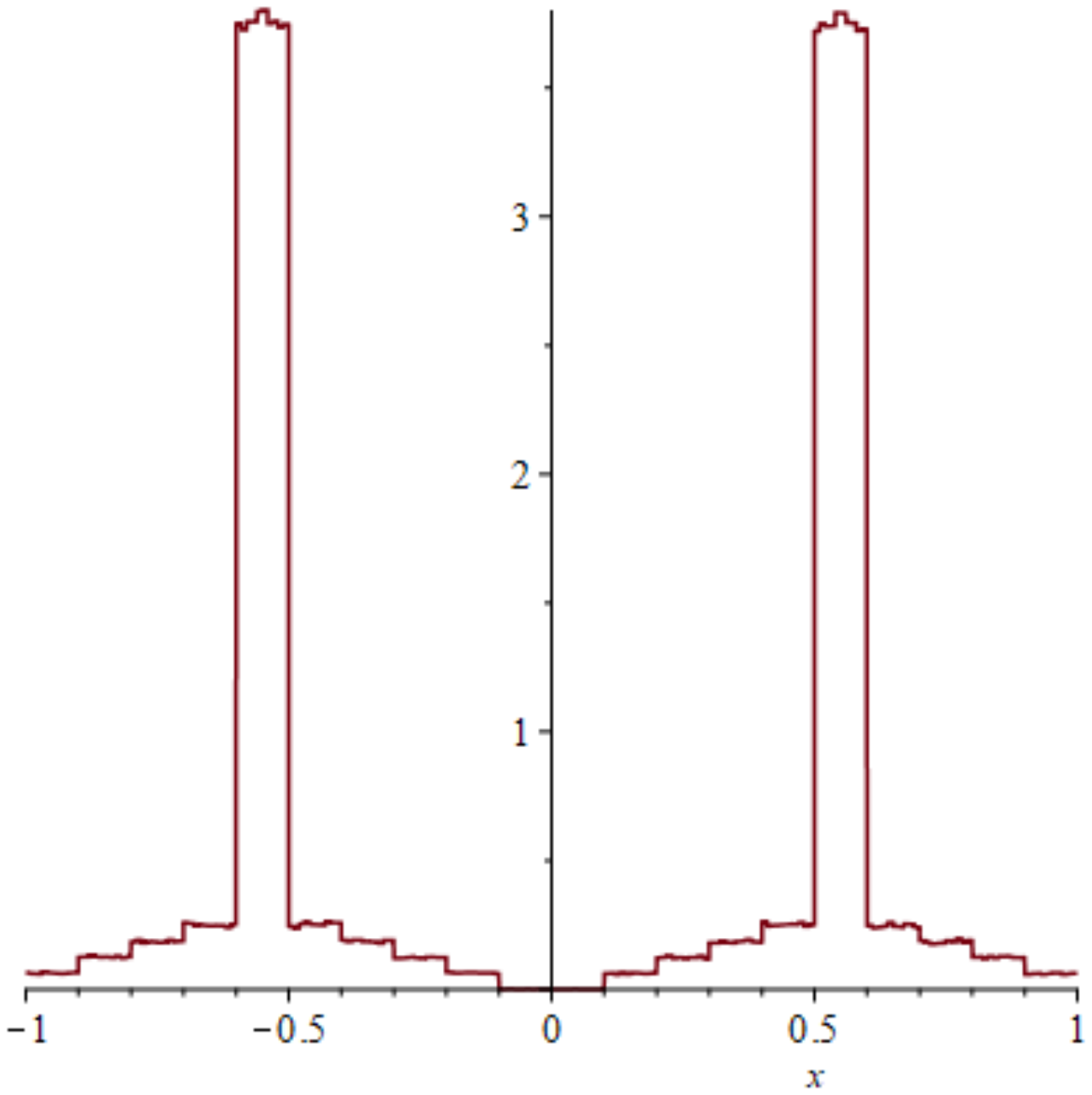}
  \caption{Metastable map $T_\eps$ preserves unique density\hskip 3 cm  $f_\eps \sim \frac 12 f_A+\frac 12 f_B$}
  \label{fig:invariant_density_metastable_10}
\end{figure}

The computer simulation of the density $f_\eps \sim \frac 12 f_A+\frac 12 f_B$ is shown in Figure
\ref{fig:invariant_density_metastable_10}.

\newpage

\section{ Conflict of Interest} On behalf of all authors, the corresponding author states that there is no conflict of interest.

\bigskip 

\end{document}